\documentclass{amsart}
\usepackage{amsmath}
\usepackage{amssymb}
\usepackage{amscd}
\usepackage{graphicx}
\usepackage{latexsym}
\usepackage{psfrag}
\usepackage{framed, color}
\newtheorem{thm}{Theorem}[section]

\theoremstyle{definition}

\theoremstyle{remark}

\newtheorem{question}[thm]{Question}

\newcommand{\M}{\textrm{Mod}}
\newcommand{\CPb}{\overline{\mathbb{CP}}{}^{2}}
\newcommand{\CP}{{\mathbb{CP}}{}^{2}}

\newcommand{\Int}{\mathop{\mathrm{Int}}\nolimits}

\newcommand{\C}{\mathbb{C}}
\newcommand{\R}{\mathbb{R}}
\newcommand{\Z}{\mathbb{Z}}

\def \x {\times}

\begin{document}

\title[Simplified broken Lefschetz fibrations and trisections] {Simplified broken Lefschetz fibrations and trisections of $4$--manifolds}
\vspace{0.2in} 

\author[R. \.{I}. Baykur]{R. \.{I}nan\c{c} Baykur}
\address{Department of Mathematics and Statistics, University of Massachusetts, Amherst, MA 01003-9305, USA}
\email{baykur@math.umass.edu}

\author[O. Saeki]{Osamu Saeki}
\address{Institute of Mathematics for Industry,
Kyushu University, Motooka 744, Nishi-ku, Fukuoka 819-0395,
Japan}
\email{saeki@imi.kyushu-u.ac.jp }

\begin{abstract}
Shapes of four dimensional spaces can be studied effectively via maps to standard surfaces. We explain, and illustrate by quintessential examples, how to simplify such generic maps on \mbox{$4$--manifolds} topologically, in order to derive simple decompositions into much better understood manifold pieces. Our methods not only allow us to produce various interesting families of examples, but also to establish a correspondence between simplified broken Lefschetz fibrations and simplified trisections of closed, oriented $4$--manifolds.
\end{abstract}

\maketitle

\setcounter{secnumdepth}{2}
\setcounter{section}{0}

% ==========================================================================================================

\section{introduction}
There is a long and rich history of studying geometry and topology of spaces by looking at maps between them. For a $4$--dimensional manifold, generic maps to surfaces allow one to foliate it by surfaces, some of which are pinched along embedded loops on them. Two of the most paramount classes of maps, which received tremendous attention in recent years, are (broken) Lefschetz fibrations and trisected Morse $2$--functions. Both yield decompositions of the ambient $4$--manifold into much simpler pieces, such as symplectic fibrations or thickened handlebodies, allowing one to bring a hefty combination of ideas and techniques from complex and symplectic geometry, classical $3$--manifold topology, and geometric group theory. 

In this article, we will approach broken Lefschetz fibrations and trisections from the vantage point of singularity theory, focusing on how to construct much simplified versions of these maps/decompositions through topological modifications of generic maps. Most of the topological beautifications, we perform by homotopies of maps, guided and argued via diagrammatic representations of their singular \linebreak images ---which we hope will make our rather combinatorial arguments accessible to a broader audience. The reader is invited to glance over some of the figures below.

Our main goal is to demonstrate, with illuminative examples, how naturally and easily such simplified maps and trisections arise on $4$--manifolds. Here we will show how to pass from a simplified broken Lefschetz fibration to a simplified trisection and back, without increasing the fibration/trisection genus much. Further, we will authenticate infinite families of examples for smallest possible genera. Background results and their complete proofs, which are of fairly technical nature, are given in our more extensive work in \cite{BSsimplified}.

\section{Maps with elementary singularities}

First, we introduce the classes of maps we are interested in. Hereon, $f\colon X \to \Sigma$ is a smooth map from a closed, connected, oriented, smooth $4$--manifold $X$ to a compact, connected, oriented surface $\Sigma$. 

\subsection*{Generic maps}
The map $f$ is said to have a \emph{fold singularity} at $y$, if there are local coordinates around $y$ and $f(y)$, with respect to which the map can be written as
\[(t, x_1, x_2, x_3) \mapsto (t, \pm x_1^2 \pm x_2^2 \pm x_3^2) ,\]
and a \emph{cusp singularity} if 
\[(t, x_1, x_2, x_3) \mapsto (t, x_1^3 + t x_1 \pm x_2^2 \pm x_3^2) . \]
A fold or a cusp point is \emph{definite} if the coefficients of all quadratic terms in the corresponding local model are of the same sign, \emph{indefinite} otherwise. 

Fold and cusp points constitute a $1$--dimensional submanifold $Z_f$ of $X$, namely, a disjoint union of finitely many arcs and circles of folds, and finitely many cusps as the end points of fold arcs.  The singular image $f(Z_f)$ is, generically, a collection of cusped immersed curves on $\Sigma$ with transverse double points along fold points. When crossing over the image of a component of $Z_f$ from one side to the other, the fibers change by an index--$i$ handle attachment, where $i$ is the number of negative coefficients of the quadratic terms in the local model chosen with compatible orientation. 

By Thom transversality \cite{Thom}, any smooth map $f\colon X \to \Sigma$ can be approximated arbitrarily well by a \emph{generic map}, which has only fold and cusp singularities. Generic maps hence played a vital role in \emph{singularity theory}, following the foundational works of Whitney, Thom and Arnold since the 1950s.

\subsection*{Broken Lefschetz fibrations}
The map $f$ is said to have a \emph{Lefschetz singularity} at a point $y \in X$, if there are orientation-preserving local coordinates around $y$ and $f(y)$, with respect to which it conforms to the complex model
\[(z_1, z_2) \mapsto z_1 \, z_2 \, .\]

Lefschetz critical points constitute a $0$--dimensional submanifold $C_f$ of $X$, to wit a disjoint union of finitely many points. A \emph{broken Lefschetz fibration} $f\colon X \to S^2$ is then a surjective map with only Lefschetz and indefinite fold singularities. Any singular fiber is obtained by collapsing embedded loops on a regular fiber; each Lefschetz singularity yields an isolated node, whereas each indefinite fold circle yields a parametrized family of locally pinched surfaces.

This class of maps were first introduced by Auroux, Donaldson and Katzarkov in \cite{ADK}, as a generalization of honest \emph{Lefschetz fibrations} (without folds), which have become central objects in symplectic and contact geometry after the pioneering works of Donaldson, Gompf, Seidel and Giroux since mid-1990s. While only symplectic $4$--manifolds admit Lefschetz fibrations (with $C_f \neq \emptyset$), on any $X$ there are generic maps that can be homotoped to broken Lefschetz fibrations.

\subsection*{Trisections}
A trisection of a $4$--manifold $X$ is a decomposition into three $4$--dimensional $1$--handlebodies (thickening of a wedge of circles) meeting pairwise in $3$--dimensional \mbox{$1$--handlebodies,} and all three intersecting along a closed, connected, orientable surface. 

Trisections were introduced by Gay and Kirby in \cite{GKtrisections} as natural analogues of Heegaard splittings of \mbox{$3$--manifolds}. Just like how Heegaard splittings correspond to certain Morse functions, which are generic maps to the $1$--dimensional disk $D^1$, trisections correspond to certain generic maps to the $2$--dimensional disk $D^2$, called \emph{trisected Morse $2$--functions}. This class of maps are characterized by the following: up to isotopy, they have a single \emph{definite} fold circle mapped to the boundary $\partial D^2$, the base $D^2$ can be non-singularly foliated by rays from a regular value (say the origin) to $\partial D^2$, each intersecting the indefinite singular image always in the direction of index--$2$ handle attachments. In addition, three of these rays split $D^2$ into three sectors, where there is at most one cusp on each singular arc image in a sector, and the total number of cusps in the sectors are equal. Furthermore, the singular arcs are situated inside. (See Figure~\ref{fig:general_vs_simplified_trisection} below for a less wordy description.) We will simply call a trisected Morse $2$--function on $X$ a \emph{trisection} of $X$, while keeping in mind that non-isotopic trisected Morse $2$--functions may give rise to equivalent trisection decompositions. For $g'$ and $k'$ the numbers of indefinite fold arcs and indefinite fold arcs without cusps in each sector, respectively, we get a so-called \emph{$(g',k')$--trisection} of $X$, where $g'$ is the \emph{genus of the trisection}.  The preimages of these three sectors are the three $4$--dimensional $1$--handlebodies of the trisection decomposition.  

Any $4$--manifold $X$ admits generic maps that can be homotoped to a trisection. Even more remarkably, like the Reidemeister--Singer theorem for Heegaard splittings of \mbox{$3$--manifolds,} trisections of $4$--manifolds are unique up to an innate \emph{stabilization} operation \cite{GKtrisections}.

\section{Maps with simplified topologies}

We now describe special subclasses of broken Lefschetz fibrations and trisections, which have simpler topologies.

\subsection*{Simplified broken Lefschetz fibrations} 
A broken Lefschetz fibration $f\colon X \to S^2$ is said to be \emph{simplified}, if it satisfies the following additional properties: $f$ is injective on $Z_f \cup C_f$, all fibers and $Z_f$ (possibly empty) are connected, and $f(C_f)$ lies on the component of $S^2 \setminus f(Z_f)$ with higher genus fibers. The \emph{genus of $f$} is that of a higher genus regular fiber.

This subclass of broken Lefschetz fibrations were introduced by the first author in \cite{BaykurPJM}. The underlying topology is simple: either we get a genus--$g$ Lefschetz fibration over $S^2$ (when $Z_f = \emptyset$), or $f$ decomposes into a genus--$g$ Lefschetz fibration over a $2$--disk, a trivial genus--$(g-1)$ bundle over a $2$--disk, and a fibered cobordism between them realized by a single \emph{round $2$--handle} ($S^1$ times a $2$--handle). There are two remarkable advantages to this simplified picture. First, it induces a simple handle decomposition of $X$; see e.g. Figure~\ref{fig:genus_one_sblfs} below. Second, it makes it possible to recast the fibration algebraically in terms of Dehn twist factorizations of mapping classes \cite{BaykurPJM}.
 
Let us digress a little on the latter aspect. Let $\M(\Sigma_g)$ denote the mapping class group of orientation-preserving diffeomorphisms of the genus--$g$ surface $\Sigma_g$, and $\M _c \, (\Sigma_g)$ denote its subgroup that consists of  mapping classes which stabilize an embedded loop $c$. Let $T_{c_i}$ denote a positive (right-handed) Dehn twist along $c_i$. Then, associated to a simplified broken Lefschetz fibration, there is an ordered tuple of cycles $(c; c_1,\ldots,c_k)$, with $c, c_i$ embedded loops on $\Sigma_g$, such that
\[ \mu= T_{c_k} \cdot \ldots \cdot T_{c_1}, \ \mu(c) = \pm c  \ , i.e. \ \, \mu \in \M _c (\Sigma_g) \,, \text{and} \]
\[ \mu \in \text{Ker}(\Phi_c: \M_c(\Sigma_{g}) \to \M(\Sigma_{g-1})) \, , \]
where the latter is the homomorphism induced by first cutting $\Sigma_g$ along a non-separating loop $c$ and then gluing disks to the two new boundary components. Conversely, any ordered tuple of curves satisfying these algebraic conditions yield a genus--$g$ simplified broken Lefschetz fibration, where $\mu$ is the global monodromy of the genus--$g$ Lefschetz fibration over $D^2$ with $k$ singular points, and $c$ is the loop surgered fiberwise to match the (trivial) monodromy of the one smaller genus surface bundle over $D^2$. Here $X$ is recaptured uniquely provided $g \geq 3$, or otherwise with some additional data for identifying the corresponding end of the round handle cobordism with the boundary of a $g=0$ or $1$ fibration over the $2$--disk (parametrized roughly by $\Z_2$ or $\Z^2$, respectively).

\subsection*{Simplified trisections}
\begin{figure}[htbp!] 
\begin{center}
\includegraphics[width=\linewidth,height=0.3\textheight,
keepaspectratio]{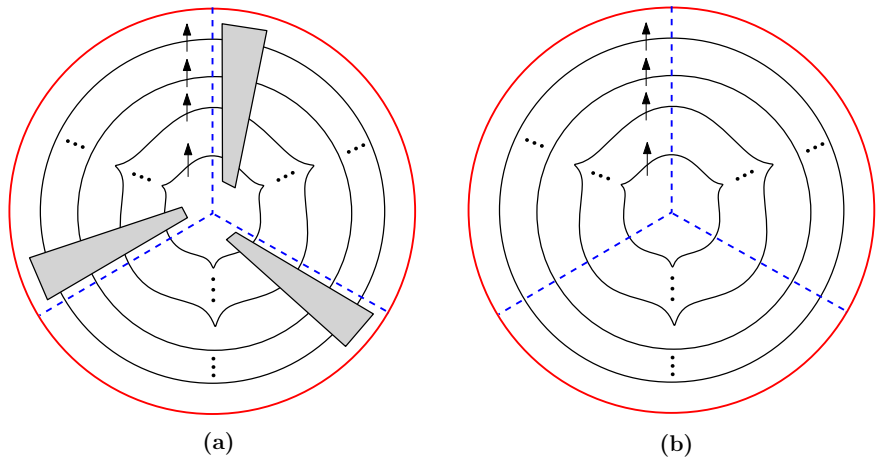}
\caption{(a) The singular image of a generic map corresponding to a trisection. The outermost circle is the definite fold, where the three gray boxes contain arbitrary Cerf graphics with intersections between folds. The three dashed-rays (given in blue) divide the base disk into three sectors, preimages of which correspond to the three $4$--dimensional $1$--handlebody pieces of the trisection. (b) Singular image of a simplified trisection with "only circles and triangles". In both pictures, dots indicate repeated patterns.}

\label{fig:general_vs_simplified_trisection}
\end{center}
\end{figure}

A trisection is said to be \emph{simplified} if the singular image of $f$ is embedded and cusps only appear in triples --like a \emph{triangle}-- in innermost fold circles. This is in great contrast with a general trisection, which has the so-called \emph{Cerf boxes} in between the three sectors of the base disk, where folds can cross each other arbitrarily (and therefore, the images of some indefinite fold circles might wind around the origin multiple times). 
Compare the singular images given in Figure~\ref{fig:general_vs_simplified_trisection}, where the arrows indicate the index--$2$ fiberwise handle attachments, and the definite fold is given in red. 

Simplified trisections were introduced recently in \cite{BSsimplified}. The main difference between a general trisection and a simplified one is manifested in what one might call the \emph{hierarchy of handle slides}. The inverse image, under an arbitrary trisection, of any radial cut of the base disk from the origin to its boundary (say, avoiding the cusps) is a genus--$g'$ handlebody, obtained by compressing $g'$ disjoint embedded loops $c_i$ on the genus--$g'$ fiber over the origin. As expected, these $c_i$ come from the fiberwise \mbox{$2$--handle} attachments prescribed by the corresponding fold arcs the ray crosses over. In general, when we move the ray across a  \emph{non-trivial} Cerf box, the corresponding $c_i$ may slide over each other in arbitrary fashion; even the roles of any two $c_i$ and $c_j$ may be swapped. Whereas for a simplified trisection, a $c_i$ slides over $c_j$ only if $i > j$. To put it loosely, these handle slides may occur only in ``upper-triangular'' fashion.

\subsection*{The existence}

With all the necessary definitions in place, we can now quote our main result from \cite{BSsimplified} motivating this work:

\begin{thm}[Existence] \label{existence} 
Given any generic map from a closed, connected, oriented, smooth $4$--manifold $X$ to $S^2$, there are explicit algorithms to modify it to a simplified broken Lefschetz fibration, as well as to a simplified trisection. Therefore, any $X$ admits simplified broken Lefschetz fibrations and simplified trisections.
\end{thm}

There is of course an abundance of generic maps from any $X$ to $S^2$. Our algorithm modifies the given generic map through various homotopies to produce a simplified broken Lefschetz fibration (over $S^2$), whereas to produce a simplified trisection (over $D^2$), we in addition apply a topological modification. We will demonstrate some of these modifications in our proofs and examples shortly. The mere existence of simplified broken Lefschetz fibrations was already known by a potpourri of arguments from handlebody theory, contact geometry, and singularity theory. (These lines of arguments do not provide explicit constructions, as they involve implicit steps such as invoking Giroux's stabilization result for contact open books.) The existence of simplified trisections is new.  

We can then derive these two types of simplified maps from one another, and the next theorem aims to do it in the most economical way for the genus of the resulting broken Lefschetz fibration or trisection, i.e. by keeping it as small as we can (but not to say one cannot do better for specific examples):

\begin{thm}[Correspondence] \label{correspondence}
If there is a genus--$g$ simplified broken Lefschetz fibration $f\colon X \to S^2$, with $k \geq 0$ Lefschetz critical points, and $\ell \in \{0, 1\}$  components of $Z_f$, then there is an associated simplified \mbox{$(g',k')$--trisection} of $X$, where
\[ g'=2g+k-\ell+2 \ \text{and} \ k'=2g-\ell \, . \] Conversely, if $X$ admits a simplified $(g', k')$--trisection, then there is an associated genus--$g$ simplified broken Lefschetz fibration \mbox{$f\colon X \to S^2$}, with $k$ Lefschetz singularities, where 
\[g=g'+3 \ \text{and} \ k=5g'-3k'+8 \, . \]
\end{thm}

%\noindent 
We will discuss the proof of this theorem below. More general versions of both directions are proved in \cite{BSsimplified}. Our proof will make use of the homotopy moves we discuss next.

\section{Homotopies of generic maps and base diagram moves}

The \emph{base diagram} of a map $f\colon X \to \Sigma$ with generic and Lefschetz type singularities, is the pair $(\Sigma, f(Z_f \cup C_f))$. We normally orient the image of any indefinite fold arc or circle by a small transverse arrow in the direction of the fiberwise index--$2$ handle attachment. We depict the definite fold circles in red, and usually it will be obvious from the rest of the diagram in which direction the fiberwise index--$3$ handle attachment is, since  one side would have empty preimage. At times we will label a region by an asterisk $(*)$ to indicate that the fibers over this region are connected. In these diagrams, we denote the Lefschetz critical points by small crosses. We assume $f$ is injective on $Z_f \cup C_f$, except possibly at fold double-points.

We will perform homotopies through a sequence of \emph{base diagram moves}, viz. local modifications of a base diagram, each one of which can \emph{always} be realized by a $1$--parameter family of smooth maps (which do not change outside of this locality). While the transition happens around one point on $\Sigma$, the bifurcation of the map may occur around one point (a \emph{mono-germ move}), or two to three points (a \emph{multi-germ move}) in $X$. It turns out that only some of the possible local changes that can occur in a base diagram during a generic homotopy are \emph{always-realizable}. Yet, the bifurcations we get through always-realizable ones will be enough to obtain the desired topology for the resulting map. 

\begin{figure}[h]
\centering
\includegraphics[width=1\linewidth]{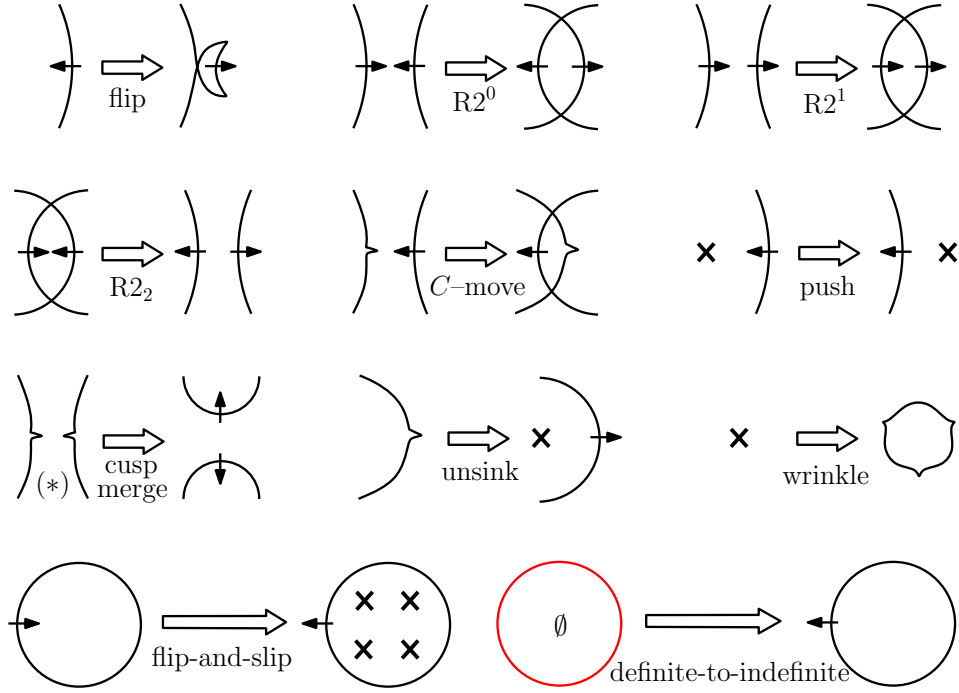}
\caption{Some always-realizable multi-germ, mono-germ and combination moves.}
\label{fig:allowable_moves}
\end{figure}

These homotopy moves have been studied in varying levels of details since mid-1960s by Levine, Hatcher--Wagoner,  Eliashberg--Mishachev, Lekili, Williams, Gay--Kirby, Behrens--Hayano, and the authors of this article.
Figure~\ref{fig:allowable_moves} lists the \emph{always-realizable moves we will employ in this article}, with the names and conventions carried on from \cite{BSsimplified}. The normal orientations for cusped arcs are always in the direction of cusps. In our arguments to follow, we will use only these always-realizable base diagram moves. (So the reader can refer to this figure as a chart of legal moves in a board game of sorts.)

The first two rows of Figure~\ref{fig:allowable_moves} consist of a mono-germ move  \emph{flip}, and multi-germ moves $\mathrm{R2}^0$, $\mathrm{R2}^1$, $\mathrm{R2}_2$, \emph{$C$--move}, and  \emph{push}. Several of these can be regarded as Reidemeister I and II type moves. There are Reidemeister III type moves as well, which play a vital role in the proof of Theorem~\ref{existence}, but are not needed for our relatively more straightforward constructions here.
The third row contains three mono-germ moves \emph{cusp merge}, \emph{unsink} and \emph{wrinkle}. Note that two cusps can be merged using any path between them in the source $4$--manifold, but here we simply use an arc with image embedded in the middle region between the two cusped arcs. When the fibers in this region are connected, one can always find such an arc. Lastly, the fourth row of Figure~\ref{fig:allowable_moves} lists two \emph{combination moves}: \emph{flip-and-slip} and \emph{definite-to-indefinite}. These involve a sequence of base diagram moves suppressed in this presentation. Importantly, almost none of these base diagram moves have always-realizable pseudo-inverses.

\section{Bridging broken Lefschetz fibrations and trisections}

Here we outline the proof of Theorem~\ref{correspondence}, in hopes to provide insight to the reader how we can use the always-realizable base diagram moves for re-arranging the underlining topology of a map to our liking. More details for the arguments below can be found in \cite{BSsimplified}.

\subsection*{From simplified broken Lefschetz fibrations to trisections}

Let $f\colon X \to S^2$ be a genus--$g$ simplified broken Lefschetz fibration with $k$ Lefschetz critical points. We will show how to derive a simplified trisection on $X$ from $f$.

\begin{figure}[htbp] 
\centering
\includegraphics[width=\linewidth,height=0.33\textheight,
keepaspectratio]{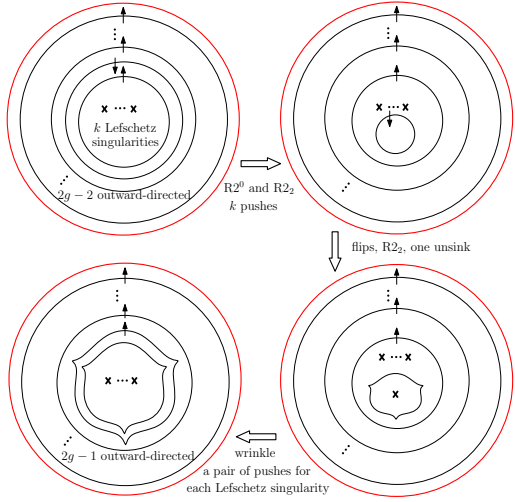} 
\caption{The singular image for the generic map $h$ obtained by folding the simplified broken Lefschetz fibration, followed by the sequence of base diagram-moves for turning it into a simplified trisection. Repeating the last step for each Lefschetz singularity, we end up with a simplified trisection.}
\label{fig:sblf_to_sts}
\end{figure}

First assume that $Z_f \neq \emptyset$. Decompose the base $S^2$ into two disks $D^2_+$ and $D^2_-$ \, such that the entire singular image $f(Z_f \cup C_f)$ lies in the \emph{interior} of $D^2_+$ and $f(Z_f)$ is parallel to the equatorial circle $E= \partial D^2_+ = - \partial D^2_-$. Identifying the base $S^2$ with the unit $2$--sphere in $\R^3$, so that $E$ maps to the boundary of the unit \mbox{$2$--disk} $D^2$ in $\R^2 \times \{0\}$, we consider the map that projects $S^2$ onto $D^2$. We ``fold'' the given simplified broken Lefschetz fibration by composing it with the described projection, in order to derive a new map to $D^2$. A careful perturbation of this map (which was constructed from scratch in \cite{BSsimplified}) is a generic map $h \colon X \to D^2$, with an \emph{embedded} singular image as follows: a definite fold along $\partial D^2$ and $2g$ boundary parallel, concentric, indefinite fold circles, enclosing Lefschetz critical points in the center. The innermost fold circle is the one that corresponds to the original $f(Z_f)$, and it is directed outwards. The next one is directed inwards, and all others outwards. At this point the fibers over any region enclosed by the inward-directed indefinite fold circle have two connected components: one coming from the preimage of a point in $D^2_+$ and the other from the preimage of the corresponding point in $D^2_-$. See Figure~\ref{fig:sblf_to_sts} for the base diagram of $h$.

We can now apply base diagram moves to turn $h$ into a trisection. Recall that we use the terminology from \cite{BSsimplified}; the names we call out for the moves can be found in Figure~\ref{fig:allowable_moves}. 

First, using an $\mathrm{R2}^0$ and an $\mathrm{R2}_2$ move, we can change the order of the innermost two circles. The innermost indefinite fold circle of the new map is directed inwards, all others outwards. Push all the Lefschetz singularities across this circle, so it now bounds a disk with no singularity inside. Apply two flips, and then an $\mathrm{R2}_2$ move to revert this circle to an outward-directed one, now with four cusps. Push back all the Lefschetz singularities into the region bounded by it. At this point, all the indefinite fold circles are directed outwards, whereas all the Lefschetz singularities are contained in the innermost region.  

What remains is to arrange the triple-cusped indefinite fold circles, as shown in Figure~\ref{fig:sblf_to_sts}. Unsinking one of its four cusps, the innermost fold circle becomes a triple-cusped one. Wrinkle one of the Lefschetz singularities to produce the next triple-cusped circle. We push all remaining Lefschetz 
singularities into the region bounded by this triple-cusped circle, and repeat the same procedure until we exhaust all the Lefschetz singularities. We end up with $k+2$ triple-cusped indefinite fold circles. The resulting map we have obtained is a simplified $(g',k')$--trisection with $g'=2g+k+1$ and $k'=2g-1$. 

If we had $Z_f= \emptyset$, the generic map $h$ we began with would not have the innermost circle that is directed outwards, but after the inward-directed circle, it would have $2g$ outward-directed fold circles instead. Then, following the same steps as above, we would get a simplified $(g',k')$--trisection with $g'=2g+k+2$ and $k'=2g$. $\square$

\subsection*{From trisections to simplified broken Lefschetz fibrations}

Let $h \colon X \to D^2$ be a simplified $(g',k')$--trisection. We will now show how to obtain a simplified broken Lefschetz fibration on $X$ from $h$. 

\begin{figure}[htbp!] 
\includegraphics[width=\linewidth,height=0.33\textheight,
keepaspectratio]{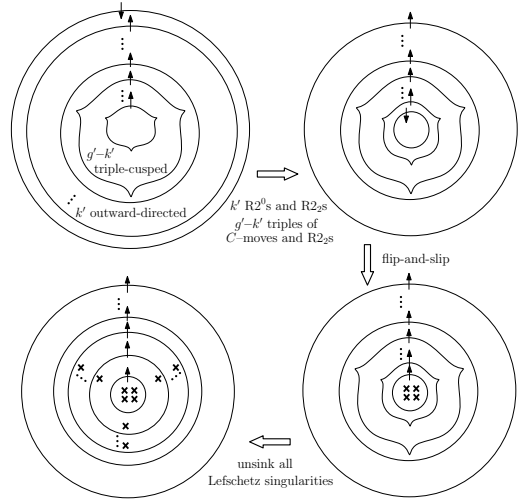} 
\caption{The singular image (on $S^2$, drawn with a point at infinity) obtained after trading the definite fold of the simplified trisection with an indefinite fold, followed by the sequence of base diagram moves for turning it into a broken Lefschetz fibration with all indefinite folds directed outwards. We then push all Lefschetz singularities into the central region, before moving onto connecting the indefinite fold locus.}
\label{fig:sts_to_sblf}
\end{figure}

Embed $D^2$ onto the northern hemisphere of $S^2$, so we view $h$ as a map to $S^2$. Applying a definite-to-indefinite move, we can trade the definite fold circle on the equator, with an indefinite one directed towards the north pole. Applying a pair of $\mathrm{R2}^0$ and $\mathrm{R2}_2$ moves repeatedly, we can move this circle across each one of the other indefinite fold circles without cusps. See Figure~\ref{fig:sts_to_sblf}. Then, by a triple of $C$--moves and a triple of $\mathrm{R2}_2$ moves, we can move it further across each one of the triple-cusped indefinite folds. Once it is the innermost circle around the north pole, we can turn it inside-out by a flip-and-slip.
Note that a regular fiber over the north pole now has genus $g'+2$. Next, we unsink \emph{all} the cusps and push the new Lefschetz singularities all the way to the innermost region around the north pole.

\begin{figure}[htbp!] 
\includegraphics[width=\linewidth,height=0.27\textheight,
keepaspectratio]{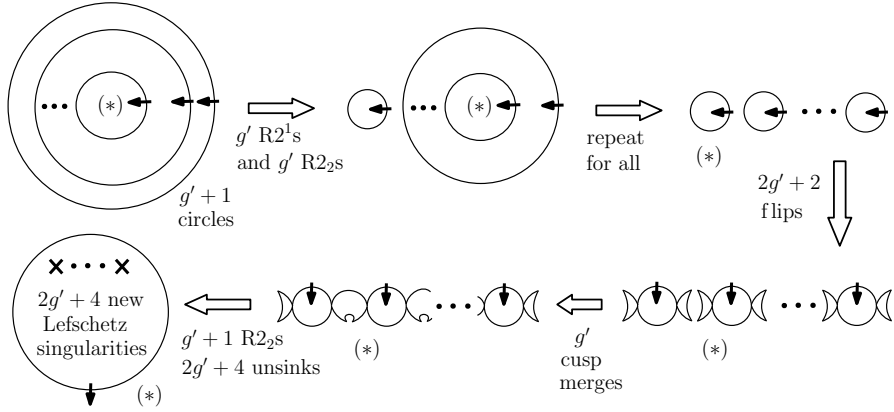} 
\caption{Merging all the indefinite fold circles into one using always-realizable base diagram moves.}
\label{fig:merging_folds}
\end{figure}

Let us view all the $g'+1$ indefinite fold circles, now none of which have cusps, in the southern hemisphere, so they are all directed inward. All the Lefschetz singularities are left in the northern hemisphere. Assume that $g' >0$. Applying $g'$  $\mathrm{R2}^1$ and $g'$ $\mathrm{R2}_2$ moves, we can push the right half of the outermost circle across all the others, so it splits from the rest. Repeating this for the others, we reach at a split collection of $g'+1$ inward-directed indefinite fold circles. See Figure~\ref{fig:merging_folds}. 

It is an easy exercise to check that throughout all the modifications we have made to $h$, every map we got so far had only connected fibers. When we have connected fibers over a region, going against the normal arrow direction of an indefinite fold, we pass to a neighboring region over which the fiber should be connected as well. This is because, tracing this path upstairs we attach an index--$1$ handle to the original regular fiber, simply increasing the genus by one. (To complete the exercise, one can for example begin with observing that after the definite-to-indefinite move, the fibers over the southern hemisphere had to be connected.) 

So if we flip each circle twice, we can merge all of them into one immersed circle after $g'$ cusp merges as in Figure~\ref{fig:merging_folds}. Applying $g'+1$ 
 $\mathrm{R2}_2$ moves we can have an \emph{embedded} indefinite fold circle directed outwards. Then, we can unsink all the cusps and push the Lefschetz singularities in the northern hemisphere across this fold circle. 
 
Since a regular fiber over the southern hemisphere is obtained by a $1$--handle attachment to a genus $g'+2$ regular fiber over the northern hemisphere (where we have not touched), its genus is $g'+3$. Along the way we created four Lefschetz singularities after flip-and-slip, $3(g'-k')$  more when we unsinked the original triples of cusps in $h$, and finally we  have $2g'+4$ more from the cusps we got in the course of merging all the circles. Hence, we have a genus--$g$ simplified broken Lefschetz fibration $f\colon X \to S^2$ with $k$ Lefschetz critical points, where $g=g'+3$ and $k= 5g'-3 k'+8$. This completes our proof. 
$\square$

\section{Small genera examples and infinite families}

We will now look at simplified broken Lefschetz fibrations and trisections of small genera. We will present some new examples for the latter, so as to classify \emph{simplified} trisections of genera at most two, and show that for each $g' \geq 3$ there are infinitely many simplified genus--$g'$ trisections.

\subsection*{Classification of small simplified broken Lefschetz fibrations} 
\mbox{Simplified} genus--$g$ broken Lefschetz fibrations of genus $g \leq 1$ are classified in \cite{BaykurKamada, Hayano1, Hayano2}, whereas a similar result seems out of reach when $g \geq 2$, even for honest Lefschetz fibrations. 

\begin{figure}[!h]
\centering
\includegraphics[width=\linewidth,height=0.3\textheight,
keepaspectratio]{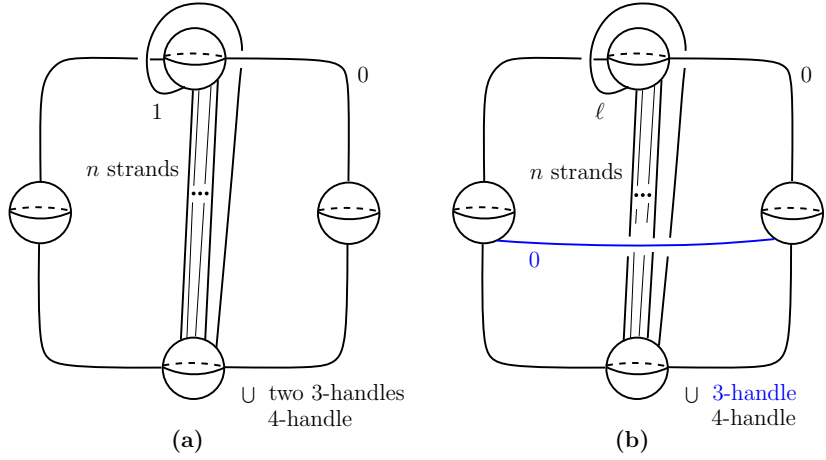}
\caption{Kirby diagrams for handle decompositions induced by genus--$1$ simplified broken Lefschetz fibrations without Lefschetz singularities. (a) $X$ is diffeomorphic to $S^1 \times L(n,1)$ for $n \geq 2$. It is $S^1 \times S^3$ for $n=1$. (b) $X$ is diffeomorphic to Pao's rational homology $4$--spheres $L_n$ and $L'_n$, $n \geq 2$, for $\ell$ even or odd, respectively. It is $S^4$ for $n=1$. Given in blue are the $2$-- and $3$--handle pair that make up the round $2$--handle of the fibered cobordism induced by the indefinite fold. (When $n=0$ and the corresponding $2$--handle is an unknot with framing $0$ instead, $X$ is diffeomorphic to $S^2 \times T^2$ in (a), and to $S^2 \times S^2 \, \# \, S^1 \times S^3$ or $X=\CP \# \CPb \, \# \, S^1 \times S^3$, for $\ell$ even or odd, respectively, in (b).)}
\label{fig:genus_one_sblfs}
\end{figure}

A genus--$0$ simplified broken Lefschetz fibration cannot have an indefinite fold, or otherwise the fibers would be disconnected. So genus--$0$ simplified broken Lefschetz fibrations are all isomorphic to holomorphic rational Lefschetz fibrations on (possibly trivial) blow-ups of complex  surfaces $S^2 \times S^2$ and $\CP \# \CPb$. A genus--$1$ simplified broken Lefschetz fibration \emph{without indefinite fold} is a classical genus--$1$ surface bundle or a Lefschetz fibration over $S^2$. If there are no critical points, these are locally trivial torus bundles on product \mbox{$4$--manifolds} $S^2 \times T^2$ (where $T^2$ denotes the $2$--torus), $S^1 \times S^3$, or $S^1 \times L(n,1)$, (where $L(n,1)$ is a Lens space), for each $n \geq 2$ \cite{BaykurKamada}. See Figure~\ref{fig:genus_one_sblfs}. If there are critical points, they are all isomorphic to holomorphic elliptic Lefschetz fibrations on (possibly trivial) blow-ups of complex elliptic surfaces $E(n)$, as shown by Kas and Moishezon in the late 1970s. An elliptic fibration on $E(n)$ has exactly $12n$ Lefschetz critical points, for each $n \geq 1$. 

The first interesting examples we get are the genus--$1$ simplified broken Lefschetz fibrations with indefinite folds. Let $a, b$ be embedded loops on $T^2$ intersecting once, and $\mu= T_a^{p} (T_a T_b)^{3q}$. It is not difficult to see that $\mu$ fixes $a$ (under the isomorphism $\M (T^2) \cong SL(2, \Z)$, the second factor is $\pm \rm{Id}$), and maps to identity under the homomorphism $\Phi_a$ discussed earlier. Thus, for each $p,q$ we get a genus--$1$ simplified broken Lefschetz fibration, and as it is, after (possibly no) blow-ups, any genus--$1$ simplified broken Lefschetz fibration becomes isomorphic to one with monodromy like this \cite{BaykurKamada}.

Surprisingly perhaps, when $\mu$ is trivial, i.e. when $p=q=0$, by varying the additional data for identifying the genus--$1$ end of the round handle cobordism with the boundary of the $g=1$ fibration, we obtain an \emph{infinite} family of examples on $4$--manifolds with distinct fundamental groups, but all with the same rational homology as the standard $4$--sphere. One can further twist this construction using the additional data for the genus--$0$ end. See Figure~\ref{fig:genus_one_sblfs}. A straightforward Kirby calculus verifies that the $4$--manifolds $L_n$ and $L'_n$  are indeed rational homology $4$--spheres, with fundamental group $\Z_n$, for $n \geq 2$. 
It turns out, $L_n$ and $L'_n$ are diffeomorphic to the rational homology $4$--spheres with effective torus actions, constructed by Pao in \cite{Pao}; see \cite{Hayano1}. The other three genus--$1$ examples with $p=q=0$, originally due to Auroux, Donaldson and Katzarkov, are on $S^4$, \mbox{$S^2 \times S^2 \, \# \, S^1 \times S^3$} and \mbox{$\CP \# \CPb \, \# \, S^1 \times S^3$}. 

The complete list of $4$--manifolds admitting genus--$1$ simplified broken Lefschetz fibration is then exhausted by (possibly trivial) blow-ups of all the $4$--manifolds we mentioned, and of $\# k\, \CP \# \CPb$ or $\# k \, S^2 \times S^2$, for any $k \geq 1$  \cite{BaykurGTM, Hayano2}.

\subsection*{Classification of small simplified trisections} 
Let us now look at the corresponding picture for \emph{simplified} trisections. General $(g',k')$--trisection decompositions are classified for $g' \leq 2$ in \cite{GKtrisections, MeierZupanGenus2}. There is only a handful of $4$--manifolds admitting them: $S^4$ for $g'=0$; $\CP$, $\CPb$ or $S^1 \times S^3$ for $g'=1$; and $S^2 \times S^2$, or connected sums of $\CP$, $\CPb$ and $S^1 \times S^3$ with two summands, for $g'=2$. We claim that the list remains the same for \emph{simplified} trisections. 

The indefinite part of the singular image of a genus--$g'$ trisection is empty when $g'=0$, and has to be either an embedded circle or an embedded triple-cusped circle when $g'=1$. So these trisections are vacuously simplified. 

As for $g'=2$, Meier and Zupan prove in \cite{MeierZupanGenus2} that \emph{any} two genus--$2$ trisections on one of the standard manifolds in the above list are the same, up to diffeomorphism. Therefore, it suffices to show that each one of these $4$--manifolds \emph{does} admit a simplified trisection. 

A few observations first. Given any map from $X$ to a surface, localizing the map over a disk with no singular image, one can always introduce a Lefschetz singularity. Furthermore, suppose we have a disk $D^2$ embedded in $X$ on which the map is a diffeomorphism onto a $2$--disk $D$ whose boundary is disjoint from the singular image. Then, we can introduce an embedded, outward-directed indefinite fold that is close and parallel to $\partial D$, and derive an extended map on $X \# \CPb$ or $X \# \, S^1 \times S^3$. The first operation is probably best known in the context of Lefschetz fibrations, where one blow-ups along a fiber. Wrinkling the Lefschetz singularity, one then gets an embedded, outward-directed, triple-cusped indefinite fold. The second one builds on the fact that $S^1 \times S^3 \setminus \Int{D^4}$ admits a map to $D^2 \cong D$ with an embedded, outward-directed indefinite fold, where a regular fiber in the center is a torus with one boundary component, and a regular fiber over $\partial D^2$ is a $2$--disk. (Hint: draw the corresponding handle diagram following \cite{BaykurPJM}, and check that the handles cancel to give $S^1 \times S^3 \setminus \Int{D^4}$.) One can then take out a fibered $D^4 \cong D^2 \times D^2$ from $X$ and extend it by the above map. We lastly note that if $X$ admits a simplified $(g',k')$--trisection, so does $\overline{X}$, the $4$--manifold with the opposite orientation. 
%(One can simply reflect the base disk.)

Using the above tricks, we can start with a genus--$1$ simplified trisection on $\CP$, $\CPb$ or $S^1 \times S^3$, and modify the $4$--manifold and the map along a generic fiber over the innermost region, or along an embedded $2$--disk as above, so as to get a genus--$2$ simplified trisection on all six $4$--manifolds that arise as the connected sums of $\CP$, $\CPb$ and $S^1 \times S^3$ with two summands. By Theorem~\ref{correspondence}, the rational fibration on $S^2 \times S^2$ (and also on $\CP \# \CPb$), regarded as simplified broken Lefschetz fibrations with \mbox{$g=\ell=k=0$,} hands us a genus--$2$ simplified trisection. This completes the proof of our claim.

\medskip
What can we say about higher genera trisections? When $g' \geq 3$, there are \emph{infinite} families of $4$--manifolds admitting genus--$g'$ simplified trisections, for fixed $g'$. For instance, the family of genus--$1$ surface bundles on $S^1 \times L(n,1)$, for $n \geq 2$, yield an infinite family of genus--$4$ simplified trisections, by Theorem~\ref{correspondence}. We will see in the next section that in fact many $L(p,q)$--bundles over $S^1$, in particular any $S^1 \times L(p,q)$, admit genus--$4$ simplified trisections.

For a sharper result, we can instead take the infinite family of \mbox{genus--$1$} simplified broken Lefschetz fibrations on rational homology \mbox{$4$--spheres} $L_n$ and $L'_n$, for $n \geq 2$, and the output of Theorem~\ref{correspondence} in this case is an infinite family of genus--$3$ simplified trisections. Let us remark that, we similarly get a genus--$3$ \emph{simplified trisection on $S^4$,} which --as a map-- is not isotopic to the standard genus--$3$ trisection used by Gay and Kirby for their stabilization result in \cite{GKtrisections}. Blow-ups of these infinite families then give infinite families of genus--$g'$ trisections for any fixed $g' \geq 3$. 

In the same fashion as our construction of genus--$2$ simplified trisections on connected sums, one can produce genus--$3$ simplified trisections on connected sums of $\CP$, $\CPb$ and $S^1 \times S^3$ with three summands, or a connected sum of either one with $S^2 \times S^2$. One can similarly get genus--$4$ examples on connected sums of lower genera trisections on these standard $4$--manifolds, and on rational homology $4$--spheres $L_n$ and $L'_n$. In addition, we have the irreducible examples on \mbox{$L(p,q)$--bundles} and $(S^1 \times S^2)$--bundles over $S^1$ we will cover in the next section, which include $S^2$--bundles over the $2$--torus $T^2$ and the Klein bottle $Kb$.

Although a complete classification of genus--$g'$ trisections seem out of reach for higher genera, it seems plausible that one can get more mileage when working with the more rigid subclass of simplified trisections, which prompts us to ask:

\begin{question} \label{q:genus3}
Which $4$--manifolds admit simplified genus--$3$ trisections? Is there any $4$--manifold, other than the ones mentioned above, which admits genus--$3$ simplified trisections? How about genus--$4$?
\end{question}

\section{More examples: from $3$--manifolds to $4$--manifolds}

Our last family of examples are on $3$--manifold bundles over the circle, and on $4$--manifolds derived from them by a standard surgery. We will show that one can easily derive a simplified broken Lefschetz fibration or a simplified trisection on these \mbox{$4$--manifolds} from any given Heegaard splitting of the $3$--manifold that is invariant under the monodromy of the bundle.

\subsection*{General constructions}
\begin{figure}[!h]
\centering
\includegraphics[width=\linewidth,height=0.3\textheight,
keepaspectratio]{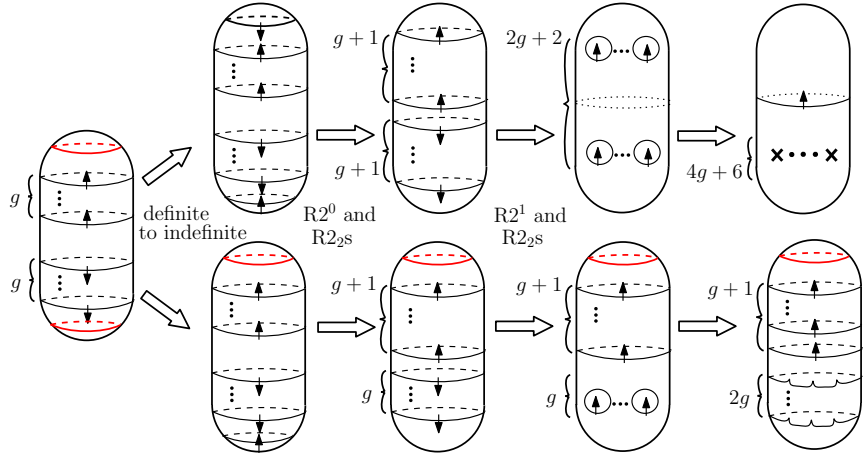}
\caption{Top: deriving a simplified broken Lefschetz fibration. In the last (longer) step, we merge all the fold circles into one and unsink the cusps that appear. Bottom: deriving a simplified trisection. In the last (even longer) step, we merge all the fold circles in the southern hemisphere into one and unsink all but three cusps. By a sequence of wrinkles and pushes applied to the Lefschetz singularities, we then turn the rest into a concentric collection of triple-cusped indefinite fold circles.}
\label{fig:heegaard}
\end{figure}

Let $Y$ be a closed, connected, oriented $3$--manifold, and let $X = S^1 \times_\varphi Y$ be the total space of a \mbox{$Y$--bundle} over $S^1$, whose monodromy is given by an orientation preserving diffeomorphism $\varphi$ of $Y$; that is
\[ X = ([-1, 1] \times Y)\, /\, (1,y) \sim (-1, \varphi(y)) \, . \]
A genus--$g$ Heegaard splitting of $Y$ gives rise to a Morse function $f_Y \colon Y\to D^1 \subset \R$ with $2g+2$ critical points, mapped injectively in a non-decreasing index order. That is, in the positive direction of $\R$, the critical values correspond to an index--$0$, then $g$ index--$1$, then $g$ index--$2$, and finally one \mbox{index--$3$} critical points. Assume that $f_Y \circ \varphi = f_Y$, so  $\varphi$ preserves the Heegaard splitting in particular. (This  is of course true for any $f_Y$ when $\varphi= \text{id}_Y$ and $X=S^1 \times Y$.) Then the product map
\[ \text{id}_{D^1} \x  f_Y \colon D^1 \x Y \to D^1 \x D^1 \, \]
descends to a generic map 
\[ \text{id}_{S^1} \x _\varphi f_Y \colon S^1 \x_\varphi Y \to S^1 \x D^1 \, . \]
Post-composing it with an embedding of the annulus \mbox{$S^1 \x D^1$} into $S^2$, we obtain a map \mbox{$f_0\colon S^1 \times_\varphi Y \to S^2$,} with singular image that consists of concentric indefinite fold circles between two definite fold circles as shown in Figure~\ref{fig:heegaard}.

Let us first derive a simplified broken Lefschetz fibration from $f_0$. By definite-to-indefinite moves we trade the two definite folds with two new indefinite folds directed towards the equator. See Figure~\ref{fig:heegaard}. Applying pairs of $\mathrm{R2}^0$ and $\mathrm{R2}_2$ moves we can move these circles so that in each hemisphere we have $g+1$ concentric circles directed towards the pole. Omitting a point in the equatorial region, we can view the whole singular image in a $2$--disk, which contains two sets of $g+1$ inward-directed concentric indefinite fold circles. Applying the first step of base diagram moves in Figure~\ref{fig:merging_folds} to each collection, we can split all as inward-directed circles. Then, following the same steps as in Figure~\ref{fig:merging_folds}, we can merge all into one, and unsink the cusps. So we have a simplified broken Lefschetz fibration $f\colon X \to S^2$. An easy book-keeping for the genera of regular fibers over the regions in each step of our modifications reveals that the genus of $f$ is $g+2$, and the number of critical points is $k= 4g+6$.

Alternatively, we can derive a simplified trisection from $f_0$. This time we only trade one of the definite folds by a definite-to-indefinite move. Applying pairs of $\mathrm{R2}^0$ and $\mathrm{R2}_2$ moves, we move the new indefinite fold to the northern hemisphere, so the singular image in the southern hemisphere now consists only of $g$ concentric inward-directed indefinite fold circles. See Figure~\ref{fig:heegaard}. Once again, we can apply the base diagram moves in Figure~\ref{fig:merging_folds} to merge all the circles in the southern hemisphere into one, but in the final step, unsink all but three of the cusps. So this new circle is a triple-cusped circle, directed towards the equator. Working with the $2g-1$ Lefschetz singularities in the same way as we did in the first part of the proof of Theorem~\ref{correspondence}, we can turn them into $2g-1$ concentric triple-cusped circles in the southern hemisphere, all directed towards the equator. The result is a simplified $(g',k')$--trisection $h\colon X \to D^2$, where $g'= 3g+1$ and $k'=g+1$.

\medskip
The constructions above are variations of those we had in \cite{BSsimplified} for the particular case of $S^1 \times Y$. Another variation, the idea of which is due to Jeff Meier, provides examples on \emph{spun $4$--manifolds}, i.e. \mbox{$4$--manifolds} obtained by surgering out $S^1 \times D^3$ from $S^1 \times Y$ and gluing in $D^2 \times S^2$ first introduced by Gordon in the 1970s. (Here there are two choices for the gluing: the end result of the gluing via the non-trivial one is usually called the \emph{twist spun}.) For this variation, instead of a definite-to-indefinite move, we remove an  $S^1 \times D^3$ neighborhood of the definite fold circle, and glue in a $D^2 \times S^2$. The map extends without any new singularity. (A similar idea was used by the second author in \cite{Saeki} originally to eliminate the definite fold from a generic map.) This results in a new $4$--manifold $X'$, which is a spun of $Y$ in this case. 

The same trick applies to any $4$--manifold $X'$ we can derive from the $3$--manifold bundle $X=S^1 \times_\varphi Y$ by a similar surgery along a $S^1 \times D^3$ neighborhood of an appropriate section of the bundle. This is possible, since the diffeomorphism $\varphi$ preserves the critical point of index $0$. The map we need to simplify here now has one less fold circle directed towards the north pole. Following the same steps as before (replacing $g+1$ with $g$ in the bottom of Figure~\ref{fig:heegaard}), we obtain a simplified $(g',k')$--trisection of $X'$, where $g'=3g$ and $k'=g$.

\subsection*{More examples of small simplified trisections}
Through the constructions above, we can produce many more simplified trisections of genus three or genus four.

Taking $Y=S^1 \times S^2$, we can obtain simplified trisections on $S^2$--bundles over $T^2$ and $Kb$. Let us explain how. When the base is $T^2$, the total spaces we get, up to diffeomorphisms, are the ruled surfaces $S^2 \times T^2$ and $S^2 \, \tilde{\times} \, T^2$. When the base is $Kb$, the orientable total spaces we get are $S^2 \times_\tau Kb$ and $S^2 \, \tilde{\times}_\tau \, Kb$. Here, $S^2 \times_\tau  Kb$ is the quotient of $S^2 \times (S^1 \times S^1)$ by the orientation preserving involution $\tau(z, a, b) =(\bar{z}, -a, \bar{b})$, where we identify $S^2$ and $S^1$ factors with $\C \cup \{\infty\}$ and the unit circle in $\C$, respectively. The twisted versions of $S^2 \times T^2$ and $S^2 \times_\tau Kb$ are then derived using the generator of \mbox{$\pi_1(SO(3)) \cong \Z_2$.} Now, how to get $S^2 \times T^2 \cong S^1 \times (S^1 \times S^2)$ is evident. For  a non-trivial example, let $\varphi_1$ be the monodromy diffeomorphism of $S^1 \times S^2$ defined by $\varphi_1(a, z) = (a, az)$. We can easily find a $\varphi_1$--invariant Morse function $f_Y$ of $Y=S^1 \times S^2$, which corresponds to a genus--$1$ Heegaard splitting. Therefore, $X = S^1 \times_{\varphi_1} (S^1 \times S^2)$, which can be seen to be diffeomorphic to the ruled surface $S^2 \, \tilde{\times} \, T^2$, admits a simplified $(4, 2)$--trisection. If instead we take the diffeomorphism $\varphi_2$ defined by $\varphi_2(a, z) = (\bar{a}, \bar{z})$, we get a simplified $(4,2)$--trisection on $X = S^1 \times_{\varphi_2} (S^1 \times S^2)$, which is diffeomorphic to $S^2 \, \times_\tau Kb$. Lastly, taking $\varphi_3 = \varphi_1 \circ \varphi_2$ yields a simplified $(4, 2)$--trisection on $X = S^2 \, \tilde{\times}_\tau\, Kb$.

Taking $Y=L(p,q)$, we obtain genus--$4$ simplified trisections on an infinite family of $4$--manifolds; the same construction applied to genus--$1$ Heegaard splittings of $L(p,q)$, for $p,q$ pairs yielding distinct $L(p,q)$ up to --possibly orientation reversing-- diffeomorphisms, generate simplified \mbox{$(4,2)$--trisections} on pairwise homotopy inequivalent product $4$--manifolds $S^1 \times L(p,q)$. As for non-trivial Lens space bundles $S^1 \times_\varphi L(p,q)$, there is essentially a unique non-trivial monodromy diffeomorphism, which up to isotopy can fix a Morse function associated to the genus--$1$ Heegaard splitting of $L(p,q)$. It can be described as follows: for $L(p, q) = (S^1 \times D^2) \cup (S^1 \times D^2)$ the standard genus--$1$ splitting, let $\varphi$ be the diffeomorphism defined by the diffeomorphism $(a, z) \mapsto (\bar{a}, \bar{z})$ on each solid torus $S^1 \times D^2$. In this case, we get a simplified $(4,2)$--trisection on $X = S^1 \times_{\varphi} L(p, q)$, which is the union of two twisted $D^2$--bundles over $Kb$. 

Together with the genus--$1$ Heegaard splitting of $Y=S^3$, these indeed exhaust all the $Y$--bundles over $S^1$ admitting genus--$4$ simplified trisections constructed in the above manner.

\smallskip
Recall that our second construction produces a simplified trisection on a $4$--manifold $X'$ that is the result of a surgery on a \mbox{$3$--manifold} bundle over $S^1$. Taking $3$--manifolds with genus--$1$ Heegaard splittings yield examples of simplified \linebreak {$(3,1)$--trisections.} In particular, genus--$1$ Heegaard splittings of Lens spaces $L(p,q)$ give rise to an infinite family of genus--$3$ simplified trisections on spun $4$--manifolds, which turn out to be the same as Pao's manifolds $L_n$, $L'_n$ for $n=p$ (see \cite{Meier} for an exposition) ---making Question~\ref{q:genus3} all the more curious! On the other hand, any $4$--manifold $X'$ derived from a non-trivial Lens space bundle $X=S^1 \times_\varphi L(p,q)$ that we have considered above by surgering out an $S^1 \times D^3$ and gluing in a $D^2 \times S^2$ turns out to be diffeomorphic to $S^4$ \cite{Teragaito}. 

\bigskip
We finish with a natural question: 

\begin{question} 
Is there any $4$--manifold which admits a trisection, but not a simplified one of the same genus? 
\end{question}

Defining the \emph{minimal trisection genus} (resp. \emph{minimal simplified trisection genus}) of a $4$--manifold $X$ as the smallest genus of a trisection (resp. simplified trisection) on $X$, one can similarly ask if there is a $4$--manifold whose trisection genus is smaller than its simplified trisection genus. The two are equal for \emph{all} the $4$--manifolds with (simplified) trisections of genus $g' \leq 4$ we have discussed in this article.

\bigskip
\noindent \textbf{Acknowledgements: } The authors thank Kenta Hayano for his careful comments on a draft of this manuscript. The first author was partially supported by the NSF grant DMS-1510395. The second author has been supported in part by JSPS KAKENHI Grant Numbers JP23244008, JP23654028, JP15K13438, JP16K13754, JP16H03936, JP17H01090, JP17H06128.

\bigskip
%\bibliographystyle{plain}
%\bibliography{references}

\end{document}